\documentclass[12pt]{amsart}
\usepackage{amssymb}
\usepackage{amsfonts}
\usepackage{amsthm}
\usepackage{amsmath}

\newtheorem{thm}{Theorem}[section]

\newcommand{\R}{\mathbb R}

\newcommand{\E}{\mathbb E}

\newcommand{\X}{\mathfrak X}

\renewcommand{\span}{\mathrm{span}}

\newcommand{\ds}{\displaystyle}

\title[]
{\bf Complete Integrability of a Nonlinear Elliptic System,
Generating Bi-Umbilical Foliated Semi-Symmetric Hypersurfaces in
$\R^4$}
\author{N. Kutev \and V. Milousheva}
\date{}

\begin{document}

\maketitle \thispagestyle{empty}

\begin{abstract}
We find explicitly all bi-umbilical foliated semi-symmetric
hypersurfaces in the four-dimensional Euclidean space.
\end{abstract}

\vskip 2mm
{\small {\bf Keywords:} foliated semi-symmetric hypersurfaces;
bi-umbilical semi-symmetric hypersurfaces; surfaces in the
3-dimensional sphere; non-linear elliptic systems.

\vskip 1mm {\bf 2000 MS Classification:} 35A07; 35J60; 53A07;
53A10}

\vskip 5mm
\section{Introduction} \label{S:Intr}

Let $(M^n,g)$  be an $n$-dimensional Riemannian manifold. We
denote by $T_pM^n$ the tangent space to $M^n$ at a point $p \in
M^n$ and by $\X M^n$ - the algebra of all vector fields on $M^n$.
The associated Levi-Civita connection of the metric $g$ is denoted
by $\nabla$, the Riemannian curvature tensor $R$ is defined by
$R(X,Y) = [\nabla_X, \nabla_Y] -  \nabla_{[X,Y]}; \,\, X,Y \in \X
M^n.$

A {\it semi-symmetric space} is a Riemannian manifold $(M^n,g)$,
whose curvature tensor $R$ satisfies the identity
$$ R(X,Y) \cdot R = 0$$
for all vector fields $X, Y \in \X M^n$. (Here $R(X,Y)$ acts as a
derivation on $R$).

According to the classification of Z. Szab\'{o} \cite{Szabo,
Szabo2} the main class of semi-symmetric spaces is the class of
all Riemannian manifolds foliated by Euclidean leaves of
codimension two. The foliated semi-symmetric spaces can be
considered also as Riemannian manifolds of conullity two
\cite{BKV}.

We recall that a Riemannian manifold $(M^n,g)$ is {\it of
conullity two}, if at every point $p \in M^n$ the tangent space
$T_pM^n$ can be decomposed in the form $T_pM^n = \Delta_0(p)
\oplus \Delta_0^{\bot}(p),$ where $\dim \Delta_0(p) = n-2$,  $\dim
\Delta_0^{\bot}(p) = 2$ and $\Delta_0(p)$ is the nullity vector
space of the curvature tensor $R_p$, i.e. $\Delta_0(p) = \{X \in
T_pM^n \,\, | \,\, R_p(X,Y)Z = 0; \,\,  Y,Z \in T_pM^n\}$. The
$(n-2)$-dimensional distribution $\Delta_0: p \longrightarrow
\Delta_0(p)$ is integrable and its integral manifolds are totally
geodesic and locally Euclidean. So, $(M^n,g)$ is foliated by
Euclidean leaves of codimension two.

The foliated semi-symmetric hypersurfaces in Euclidean space
$\E^{n+1}$ are studied in \cite{GanMil-IM} with respect to their
second fundamental form. They can be considered as hypersurfaces
of type number two, i.e. hypersurfaces whose rank of the second
fundamental form is equal to two everywhere. Each foliated
semi-symmetric hypersurface $M^n$ in $\E^{n+1}$ is characterized
by a second fundamental form $h = \nu_1 \, \eta_1 \otimes \eta_1 +
\nu_2 \, \eta_2 \otimes \eta_2$, where $\eta_1$ and $\eta_2$ are
unit one-forms; $\nu_1$ and $\nu_2$ are functions on $M^n$, $\nu_1
\nu_2 \neq 0$. The Euclidean leaves of the foliation are the
integral submanifolds of the distribution $\Delta_0$, determined
by the one-forms $\eta_1$ and $\eta_2$, i.e. $\Delta_0(p) = \{X
\in T_pM^n \,\, | \,\, \eta_1(X) = 0, \,\, \eta_2(X) = 0\}, \,\, p
\in M^n$.

Let $\Delta_0^{\bot}$ be the two-dimensional geometric
distribution, which is orthogonal to the distribution $\Delta_0$,
i.e. to the Euclidean leaves of the foliation of $M^n$. In case of
$\nu_1 = \nu_2$, at each point $p \in M^n$ the shape operator $A$
of $M^n$ has two equal nonzero eigenvalues corresponding to
$\Delta_0^{\bot}$, and an eigenvalue $\nu = 0$ with multiplicity
$n-2$, corresponding to $\Delta_0$. That is why the foliated
semi-symmetric hypersurfaces satisfying the condition $\nu_1 =
\nu_2$ are called \textit{bi-umbilical}.

The foliated semi-symmetric hypersurfaces in $\E^{n+1}$ are
characterized in \cite{GanMil-MathB} by the following

\begin{thm}\label{T:envelope}
A hypersurface $M^n$ in Euclidean space $\E^{n+1}$ is locally a
foliated semi-symmetric hypersurface if and only if it is the
envelope of a two-parameter family of hyperplanes in $\E^{n+1}$.
\end{thm}

Using the characterization of a foliated semi-symmetric
hypersurface as the envelope of a two-parameter family of
hyperplanes, each such hypersurface is determined by a pair of a
unit vector-valued function $l(u,v)$ and a scalar function
$r(u,v)$, defined in a domain ${\mathcal D} \subset \R^2$.

Since the vector fields $l_u$ and $l_v$ are linearly independent,
then the vector-valued function $l(u,v)$ determines a
two-dimensional surface $M^2: l = l(u,v)$, $(u,v)\in {\mathcal D}$
in $\E^{n+1}$. We use the standard denotations $E(u,v) = g(l_u,
l_u); \, F(u,v) = g(l_u, l_v); \, G(u,v) = g(l_v, l_v)$ for the
coefficients of the first fundamental form of $M^2$. According to
\cite{Chern} for each point of a regular surface there exists a
neighbourhood, in which isothermal parameters can be introduced,
i.e. $E(u,v) = G(u,v); \,\, F(u,v) = 0$. Since our considerations
are local, without loss of generality we assume that the surface
$M^2$ is parameterized locally by isothermal parameters. Then the
generated foliated semi-symmetric hypersurface $M^n$ is given by
\cite{GanMil-CR}
$$X(u, v, w^{\alpha}) = r\,l + \displaystyle{\frac{r_u}{E}\,\,l_u} + \displaystyle{\frac{r_v}{E}\,\,l_v} + w^{\alpha}\,b_{\alpha}, \quad \alpha = 1, \dots, n-2, \leqno{(1.1)}$$
where $(u,v) \in {\mathcal D}, \,\, w^{\alpha} \in \R, \,\, \alpha
= 1, \dots, n-2$, and $b_1(u,v), \dots, b_{n-2}(u,v), \,(u,v)\in
{\mathcal D}$ are $n-2$ mutually orthogonal unit vectors,
orthogonal to $\span \{l, l_u, l_v\}$.

\vskip 2mm The bi-umbilical foliated semi-symmetric hypersurfaces
 are characterized analytically \cite {GanMil-CR} by the following

\begin{thm}\label{T:basic}
Let $M^n$ be a hypersurface in $\E^{n+1}$ which is the envelope of
a two-parameter family of hyperplanes, determined by a unit
vector-valued function $l(u,v)$, represented by isothermal
parameters, and a scalar function $r(u,v)$. Then $M^n$ is
bi-umbilical if and only if \, $l(u, v)$ and $r(u,v)$ satisfy the
equalities
$$\begin{array}{l}
\vspace{2mm}
  l_{uu} - l_{vv} - \displaystyle{\frac{E_u}{E}\,l_u + \frac{E_v}{E}\,l_v}= 0;\\
\vspace{2mm}
2 l_{uv} - \displaystyle{\frac{E_v}{E}\,l_u - \frac{E_u}{E}\,l_v}= 0;\\
\vspace{2mm}
  r_{uu} - r_{vv}  - \displaystyle{\frac{E_u}{E}\,r_u + \frac{E_v}{E}\,r_v}= 0;\\
\vspace{2mm}
  2 r_{uv} - \displaystyle{\frac{E_v}{E}\,r_u - \frac{E_u}{E}\,r_v}= 0.
\end{array}$$
\end{thm}

So, the bi-umbilical foliated semi-symmetric hypersurfaces are
determined by the following system of non-linear equations for the
vector-valued function $l(u,v)$:
$$\begin{array}{l}
\vspace{2mm}
  l_{uu} - l_{vv} - \displaystyle{\frac{E_u}{E}\,l_u + \frac{E_v}{E}\,l_v}= 0;\\
\vspace{2mm} 2 l_{uv} - \displaystyle{\frac{E_v}{E}\,l_u -
\frac{E_u}{E}\,l_v}= 0,
\end{array}  \leqno{(1.2)}$$
and $l(u,v)$ satisfies the additional conditions
$$g(l_u, l_u) = g(l_v, l_v); \qquad g(l_u, l_v) = 0; \qquad g(l, l) = 1. \leqno{(1.3)}$$

By means of the identities (1.3) the system (1.2) can be rewritten
in the following way:
$$A \, l_{uu} + B\, l_{vv} = 0, \leqno{(1.4)}$$
where $\ds{A = \{I - 2\frac{l_u \otimes l_u}{l_u^2}\}, \,\, B =\{I
- 2\frac{l_v \otimes l_v}{l_v^2}\}}$. Since $\det A = - 1$, system
(1.4) has the normal form
$$l_{uu} = - A^{-1} B\,l_{vv}. \leqno{(1.5)}$$

The characteristic matrix $\det (\lambda^2 I + A^{-1}B)$ of (1.5)
has no real roots.  Indeed, calculating
$$\ds{\{\lambda^2 A + B\} = \{(\lambda^2 + 1) I + \frac{(-\lambda l_u + l_v)\otimes (-\lambda l_u + l_v)}{l_u^2}\}};$$
$$\det (\lambda^2 I + A^{-1}B) = \det A^{-1} \det (\lambda^2 I + B) = - \det (\lambda^2 A + B) = (\lambda^2 + 1)^4,$$
we get that the roots of the  characteristic matrix are $\lambda =
\pm i$.
 Hence, according to the
classification of the general systems \cite{Petr}, system (1.5) is
a non-linear elliptic system.

Let us note that (1.2) is an elliptic system in non-divergent form
and the general theory of calculus of variation can not be applied
to (1.2). That is why for solving this system we will use a
different method, which is based rather on the differential
geometry of surfaces in $\E^4$ than on the PDE methods. In such
way we find all solutions of (1.2) satisfying (1.3) in the
four-dimensional Euclidean space $\E^4$.

In Section 2 we consider the solutions of (1.2), (1.3) as
two-dimensional surfaces lying on the unit sphere $S^3(1)$ in
$\E^4$. Using the derivative formulas of these surfaces, in
Theorem 2.1 we  prove that each solution $l = l(u,v)$ of system
(1.2), (1.3) is part of  a sphere $S^2$ in  a constant hyperplane
$\R^3$ of $\E^4$.

In Section 3 we apply the result of Theorem 2.1 for giving
explicitly all bi-umbilical foliated semi-symmetric hypersurfaces
in $\E^4$.

\section{Solvability of the non-linear elliptic system}
In this section we shall consider the system
$$\begin{array}{l}
\vspace{2mm}
  l_{uu} - l_{vv} - \displaystyle{\frac{E_u}{E}\,l_u + \frac{E_v}{E}\,l_v}= 0\\
\vspace{2mm}
2 l_{uv} - \displaystyle{\frac{E_v}{E}\,l_u - \frac{E_u}{E}\,l_v}= 0\\
\vspace{2mm}
g(l_u, l_u) = g(l_v, l_v) = E\\
\vspace{2mm}
g(l_u, l_v) = 0\\
\vspace{2mm} g(l, l) = 1
\end{array} \leqno(2.1)$$
for the vector-valued function $l(u,v) = \left(l^1(u,v), l^2(u,v),
l^3(u,v), l^4(u,v)\right)$ in the 4-dimensional Euclidean space
$\E^4$.

Let $l = l(u,v)$ be a solution of (2.1), defined in a domain
${\mathcal D} \subset \R^2$. We consider the 2-dimensional surface
$M^2: l = l(u,v), \, \, (u,v) \in {\mathcal D}$ (${\mathcal D}
\subset \R^2$) in $\E^4$. $M^2$ is a surface lying on the unit
sphere $S^3(1)$ in $\E^4$, and the parameters $(u,v)$ are
isothermal ones.

\begin{thm} \label{T:umbilical}
Each solution $M^2: l = l(u,v)$ of system $(2.1)$ lies on a sphere
$S^2$ in  a constant hyperplane $\R^3$ of  $\E^4$.
\end{thm}

\vskip 2mm \noindent \textit{Proof}. Let $M^2: l = l(u,v)$ be a
solution of (2.1). The tangent space to $M^2$ at an arbitrary
point $p=l(u,v)$ of $M^2$ is $T_pM^2 = {\rm span} \{l_u, l_v\}$.
Since the vector fields $l, l_u, l_v$ are orthogonal, there exists
a unique (up to a sign) unit vector field $n(u,v)$, such that
$\{l_u, l_v, l, n\}$ form an orthogonal frame field in $\E^4$.
Using (2.1) we obtain the following derivative formulas of $M^2$:
$$\begin{array}{l}
\vspace{2mm}
 l_{uu} = \displaystyle{\frac{E_u}{2E}\,l_u - \frac{E_v}{2E}\,l_v - E\, l + c\, n};\\
\vspace{2mm}
l_{uv} = \displaystyle{\frac{E_v}{2E}\,l_u + \frac{E_u}{2E}\,l_v};\\
\vspace{2mm}
 l_{vv} = \displaystyle{- \frac{E_u}{2E}\,l_u + \frac{E_v}{2E}\,l_v - E\, l + c\, n},
\end{array} \leqno(2.2)$$
where $c(u,v) = g(l_{uu}, n) = g(l_{vv}, n)$. The equalities (2.2)
imply
$$n_u = -\displaystyle{\frac{c}{E}\,l_u}; \qquad n_v = -\displaystyle{\frac{c}{E}\,l_v}.$$
Since $n_{uv} = n_{vu}$, then $\ds{\left(\frac{c}{E}\right)_u =
\left(\frac{c}{E}\right)_v = 0,}$ and hence
$\displaystyle{\frac{c}{E}} = c_0 = const$, i.e. $c(u,v) = c_0
E(u,v)$.

Let $\{x = \ds{\frac{l_u}{\sqrt{E}}}, \,\, y =
\ds{\frac{l_v}{\sqrt{E}}}\}$ be an orthonormal tangent frame field
of $M^2$. From (2.2) we get
$$\begin{array}{ll}
\vspace{1mm}
 \nabla'_x l = x; \qquad & \nabla'_x n = - \ds{\frac{c}{E}}x;\\
\vspace{1mm}
 \nabla'_y l = y; \qquad & \nabla'_y n = - \ds{\frac{c}{E}}y.
\end{array} \leqno(2.3)$$
From (2.2) it follows that the Riemann curvature $K$ of $M^2$ is
expressed as $K = 1 + \ds{\frac{c^2}{E^2}} = 1+ c_0^2.$ Hence, the
surface $M^2$ is of constant Riemann curvature $K$.

In case of $c_0 = 0$ the normal vector field $n$ is constant and
$M^2$ lies in the constant 3-dimensional subspace $\R^3 = \span
\{l_u,l_v,l\}$ of $\E^4$. Moreover, $M^2$ lies on a sphere $S^2(1)
= S^3(1)\bigcap \R^3$.

In case of $c_0 \neq 0$ we consider the following normal vector
fields
$$n_1 = \ds{\frac{1}{\sqrt{c^2+E^2}}(c\,l + E\,n)}; \qquad n_2 = \ds{\frac{1}{\sqrt{c^2+E^2}}(-E\,l + c\,n)}. \leqno{(2.4)}$$
From (2.3) and (2.4) we get
$$\begin{array}{ll}
\vspace{2mm}
 \nabla'_x n_1 = 0; \qquad & \nabla'_x n_2 = - \ds{\sqrt{1+c_0^2}}\,x;\\
\vspace{2mm}
 \nabla'_y n_1 = 0; \qquad & \nabla'_y n_2 = - \ds{\sqrt{1+c_0^2}}\,y.
\end{array}$$
Hence, $M^2$ lies on a sphere $S^2(\frac{1}{\sqrt{1+c_0^2}})$   in
the constant 3-dimensional subspace $\R^3 = \span \{l_u,l_v,
n_2\}$ of $\E^4$. \qed

\section{Application to the geometric problem} \label{S:Bi-umbilical}

In Section 2 we proved that each solution $l = l(u,v)$ of (2.1) is
a 2-dimensional surface $M^2: l = l(u,v)$ lying on $S^3(1) \bigcap
\R^3$ for some constant hyperplane $\R^3$. The standard
parametrization of the unit sphere $S^3(1)$ in $\E^4$ is
$$S^3: z(u,v,\alpha) = \cos \alpha \, \cos u \, \cos v \, e_1 + \cos \alpha \, \cos u \, \sin v \, e_2 +
\cos \alpha \, \sin u \, e_3 + \sin \alpha\, e_4,$$ where $\{e_1,
e_2, e_3, e_4\}$ is the standard basis in $\E^4$, and $u \in [0;
2\pi), \, v \in [0; 2\pi), \, \alpha \in [0; 2\pi)$.

Let $\R^3$ be the hyperplane in $\E^4$, defined by $z^4 = b =
const$. Then the 2-dimensional sphere $S^2 = S^3(1) \bigcap \R^3$
is parameterized by
$$S^2: l(u,v) = a \cos u \, \cos v \, e_1 + a \, \cos u \, \sin v \, e_2 +
a \, \sin u \, e_3 + b\, e_4, \leqno{(3.1)}$$ where $\cos \alpha =
a, \, \sin \alpha = b$. $S^2$ is a sphere with Gauss curvature $K
= 1 + \ds{\frac{b^2}{a^2}}$, and radius $r =
\ds{\frac{1}{\sqrt{K}}} = a$. From (3.1) we get
$$E = g(l_u, l_u) = a^2; \qquad F = g(l_u, l_v) = 0; \qquad G = g(l_v, l_v) = a^2 \cos^2 u.$$

Obviously the parameters $(u,v)$ of $S^2$ are not isothermal.
After the following change of the parameters:
$$u = \ds{- \frac{\pi}{2} + 2\arctan e^x}; \qquad  v = y,$$
we get
$$\sin u = \ds{\frac{e^{2x} -1}{e^{2x} + 1}} = \tanh x; \qquad \cos u =\ds{\frac{2e^{x}}{e^{2x} + 1} = \frac{1}{\cosh x}}.$$
So, we obtain the following parametrization of $S^2$:
$$S^2: l(x,y) = \ds{\frac{1}{\cosh x} \left(a \cos y \, e_1 + a \, \sin y \, e_2 +
a \, \sinh x \, e_3 + b \cosh x\, e_4\right).} \leqno{(3.2)}$$
From (3.2) we get $E = G = \ds{\frac{a^2}{\cosh^2 x}; \,\, F =
0},$ i.e. the parameters $(x,y)$ are isothermal. In terms of the
parameters $(x,y)$  system (2.1) is rewritten in the form:
$$\begin{array}{l}
\vspace{1mm}
l_{xx} - l_{yy} = - 2 \tanh x \, l_x;\\
\vspace{1mm} 2 \,l_{xy} = - 2 \tanh x \, l_y.
\end{array} \leqno{(3.3)}$$

\vskip 2mm Now the equalities for the scalar function $r(x,y)$
have the form
$$\begin{array}{l}
\vspace{1mm}
r_{xx} - r_{yy} = - 2 \tanh x \, r_x;\\
\vspace{1mm} 2\,r_{xy} = - 2 \tanh x \, r_y.
\end{array} \leqno{(3.4)}$$

Each scalar function $r(x,y)$, satisfying (3.4), together with the
sphere $S^2: l =l(x,y)$, defined by (3.2), generate a bi-umbilical
foliated semi-symmetric hypersurface $M^3$ in $\E^4$ according to
the geometric construction, given by (1.1).

Further we are going to find all solutions of (3.4). From the
second equality of (3.4) we get $(\cosh x \,r_y)_x = 0$, which
implies $\cosh x \,r_y = g(y)$ for some smooth function $g(y)$.
Hence for fixed $x$ we get
$$r(x,y) = \ds{\frac{1}{\cosh x} \int_o^y g(t) dt + f(x)}, \leqno{(3.5)}$$
where $f(x)$ is a smooth function. The first equality of (3.4)
implies
$$g'(y) + \int_o^y g(t) dt = \cosh x\,f''(x) + 2\sinh x\, f'(x).$$
Hence
$$g'(y) + \int_o^y g(t) dt = c_4;$$
$$\cosh x\,f''(x) + 2\sinh x\, f'(x) = c_4;\leqno{(3.6)}$$
where $c_4 = const$. All solutions for $\int_o^y g(t) dt $ are
given by
$$\int_o^y g(t) dt = c_1 \cos y + c_2 \sin y + c_4, \qquad c_1, c_2 - const. \leqno{(3.7)}$$
Equality (3.6) is equivalent to the identity $(\cosh^2 x\,f'(x))'
= c_4 (\sinh x)'$, which implies
$$f(x) = c_3 \tanh x - \ds{\frac{c_4}{\cosh x}}+ c_0. \leqno{(3.8)}$$
Now from equalities (3.5), (3.7) and (3.8) we get that all
solutions of (3.4) are given by:
$$r(x,y) = c_0 + \ds{\frac{1}{\cosh x}} (c_1 \cos y + c_2 \sin y + c_3 \sinh x), \leqno{(3.9)}$$
where $c_0, \, c_1, \, c_2, \, c_3$ are constants.

Consequently, all bi-umbilical foliated semi-symmetric
hypersurfaces in $\E^4$ are given by the following formula (up to
parametrization):
$$X(x,y,w) = \ds{r\,l + \frac{\cosh^2 x}{a^2}\,r_x\,l_x + \frac{\cosh^2 x}{a^2}\, r_y \,l_y + w\,n.}  \leqno{(3.10)}$$
Calculating $r_x$ and $r_y$ from (3.9), and using (3.10) we obtain
the coordinate functions of the bi-umbilical foliated
semi-symmetric hypersurface $M^3: X = X(x,y,w)$:
$$\begin{array}{l}
\vspace{2mm}
X^1(x,y,w) = \ds{\frac{\cos y }{\cosh x} \left[ ac_0 + b\, w - \frac{b^2}{a \cosh x} (c_1 \cos y + c_2 \sin y + c_3 \sinh x)\right] + \frac{c_1}{a}};\\
\vspace{2mm}
X^2(x,y,w) = \ds{\frac{\sin y}{\cosh x}\left[ac_0 + b\, w - \frac{b^2}{a \cosh x} (c_1 \cos y + c_2 \sin y + c_3 \sinh x)\right] + \frac{c_2}{a}};\\
\vspace{2mm}
X^3(x,y,w) = \ds{ \frac{\sinh x}{\cosh x}\left[ac_0 + b\, w - \frac{b^2}{a \cosh x} (c_1 \cos y + c_2 \sin y + c_3 \sinh x)\right] + \frac{c_3}{a}};\\
\vspace{2mm} X^4(x,y,w) = -\ds{\frac{a}{b} \left[ac_0 + b\, w -
\frac{b^2}{a \cosh x} (c_1 \cos y + c_2 \sin y + c_3 \sinh
x)\right]+ \frac{c_0}{b}}.
\end{array}$$

Let us denote $f(x,y,w) =\ds{ \frac{a}{b} \,c_0 + w - \frac{b}{a
\cosh x} (c_1 \cos y + c_2 \sin y + c_3 \sinh x)}$. We consider
the vector-valued function
$$n(x,y) =  \ds{\frac{1}{\cosh x} \left(b \cos y \, e_1 + b
\, \sin y \, e_2 + b \, \sinh x \, e_3 - a \cosh x\,
e_4\right)}.$$

Then each bi-umbilical foliated semi-symmetric hypersurface $M^3:
X = X(x,y,w)$ is parameterized as follows:
$$ X(x,y,w) = f(x,y,w)\,n(x,y) + C,$$
where $C = \ds{\frac{c_1}{a}\, e_1 + \frac{c_2}{a}\, e_2 +
\frac{c_3}{a}\, e_3 + \frac{c_0}{b}\, e_4}$ is a constant vector
in $\E^4$.

\vskip 2mm
Let us note that in \cite{Gan} G. Ganchev studied the integral surfaces of the distribution $\Delta_0^{\bot}$
orthogonal to the Euclidean leaves of the foliation and proved that the integral surfaces of $\Delta_0^{\bot}$
of any bi-umbilical semi-symmetric hypersurface lie on two-dimensional spheres. Conversely, any two-dimensional sphere
generates a family of bi-umbilical semi-symmetric hypersurfaces.
Here we obtained explicitly all bi-umbilical semi-symmetric hypersurfaces using the approach of studying
these hypersurfaces as the envelopes of two-parameter families of hyperplanes and solving the corresponding
non-linear elliptic system. This geometric method of solving systems of PDE can also be applied
for the investigation of other non-linear systems arising in geometric problems (see \cite{K-M-CR}).

\vskip 4mm \textbf{Acknowledgements:} The second author is
partially supported by "L. Karavelov" Civil Engineering Higher
School, Sofia, Bulgaria under Contract No 10/2009.

\vskip 1mm \noindent
Institute of Mathematics and Informatics,\\ Bulgarian Academy of Sciences,\\
Acad. G. Bonchev Str. bl. 8, 1113, Sofia, Bulgaria\\

\noindent
e-mail: \textsf{kutev@math.bas.bg}\\
e-mail: \textsf{vmil@math.bas.bg}

\end{document}